\documentclass[a4paper, 11pt]{article}
\usepackage{amsmath, amssymb, amsthm, amsfonts}
\usepackage[dvips]{graphicx}
\newtheorem{thm}{Theorem}[section]
\newtheorem{dfn}{Definition}[section]

\begin{document}
\title{Eigenvalue Statistics of One-Face Maps}   
\author{E. M. McNicholas\\
  Department of Mathematics\\
  Willamette University\\
  900 State St\\
  Salem, Oregon 97301\\
  \texttt{emcnicho@willamette.edu}}
\date{\today}
\maketitle

\begin{abstract}
We examine the adjacency matrices of three-regular graphs
 representing one-face maps.  Numerical studies reveal that the
 limiting eigenvalue statistics of these matrices are the same as
 those of much larger, and more widely studied classes from Random Matrix Theory.  We present an algorithm
 for generating matrices corresponding to maps of genus zero, and find the eigenvalue statistics
 in the genus zero case differ strikingly from those of
 higher genus.  These results lead us to conjecture that the eigenvalue
 statistics depend on the rigidity of the underlying map, and the distribution of
 scaled eigenvalue spacings shifts from that of the Gaussian Orthogonal Ensemble to the exponential distribution
 as the map size increases relative to the genus.


\end{abstract}

\section{Background}
\label{Background}

In this paper we examine a special class of three-regular graphs. A
graph $G=\{V,E\}$ consists of a set of vertices $V$ and a set of
edges $E$, with each edge connecting two vertices.  A graph is
$k$-regular if the number of edges incident to each vertex is
exactly $k$.  The structure of a graph is encoded in the associated
adjacency matrix $A=[a_{ij}]$, where $a_{ij}$ is the number of edges
between vertex $i$ and vertex $j$. When we refer to the spectrum of
a graph, we are referring to the eigenvalues of its adjacency
matrix.

Our special class of three-regular graphs have $2N$ vertices and can
be decomposed into a cycle graph and a totally disjoint, one-regular
"toothpick" graph. The cycle graph on $2N$ vertices is the
two-regular graph in which vertex $i$ is adjacent to vertices $i+1$
and $i-1 \mod 2N$. In the standard toothpick graph, vertex $1$ is
adjacent to vertex $2$, vertex $3$ is adjacent to vertex $4$, and so
on.  The adjacency matrices of our graphs can be expressed as
$C+P^TTP$, where $C$ is the adjacency matrix of the cycle graph, $T$
is the adjacency matrix of the standard toothpick graph, and $P$ is
a permutation matrix.  Three-regular graphs of this form can be used
to represent one-face maps.

\begin{dfn}
A map (or dessin d'enfant) is an embedding of a connected, labeled
graph $X$, into a compact oriented surface $S$, with the following
properties:
\begin{enumerate}
\item The edges of the graph do not intersect (except at vertex points).
\item The complement of $X$ in $S$ is a disjoint union of open cells which are
all homeomorphic to the disk. These cells are called the {\it faces}
of the map.
\end{enumerate}
The {\it genus of a map} is the genus of the underlying surface, not
the graph-theoretic genus of the embedded graph $X$.
\end{dfn}

\begin{figure}[ht]
\begin{center}
\includegraphics{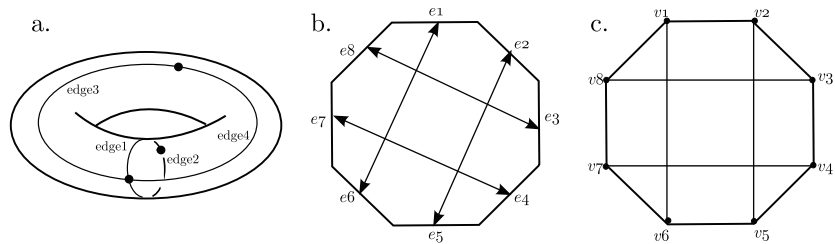}
\caption[A one-face map, its $2N$-gon representation, and the
associated three-regular graph]{a. The genus one, one-face map  b.
The corresponding $2N$-gon with edge identifications  c. The
associated three-regular graph of the form
$C+P^TTP$\label{graphrepCh1}}
\end{center}
\end{figure}

Our three-regular graphs describe a glueing of the $2N$-gon's edges.
If the toothpick portion of the three-regular graph contains an edge
connecting vertex $i$ and vertex $j$, then edge $i$ and edge $j$ are
glued together in the corresponding $2N$-gon. This glueing maps the
$2N$-gon to a graph embedded in a Riemann surface of genus $g\ge 0$.
Since the complement of this graph is homeomorphic to the unit disc,
we say the map has ``one face''. Figure \ref{graphrepCh1} shows a
one-face map, its representation by a $2N$-gon with edge
identifications, and the associated three-regular graph. All
one-face maps can be represented by a three-regular graph of the
form $C+P^TTP$.

\section{Numerical studies of one-face maps}

\begin{figure}[ht]
\begin{center}
\includegraphics{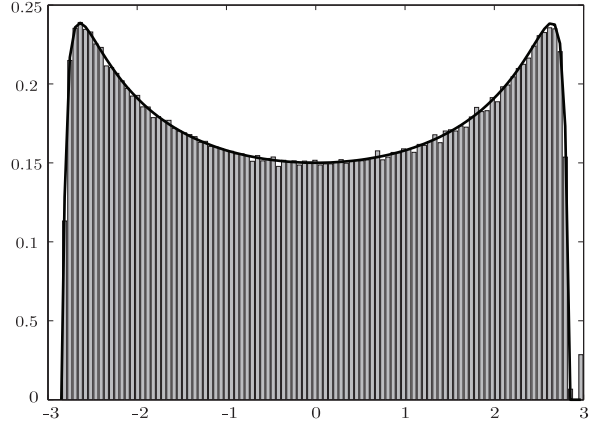}
\caption[Empirical eigenvalue density for a random sample of
matrices having the form $C+P^TTP$]{Empirical eigenvalue density for
a random sample of 300 matrices having the form $C+P^TTP$, plotted
with the McKay density function.  The corresponding graphs have 600
vertices each.\label{EmpMcKay}}
\end{center}
\end{figure}
\begin{figure}[ht]
\begin{center}
\includegraphics{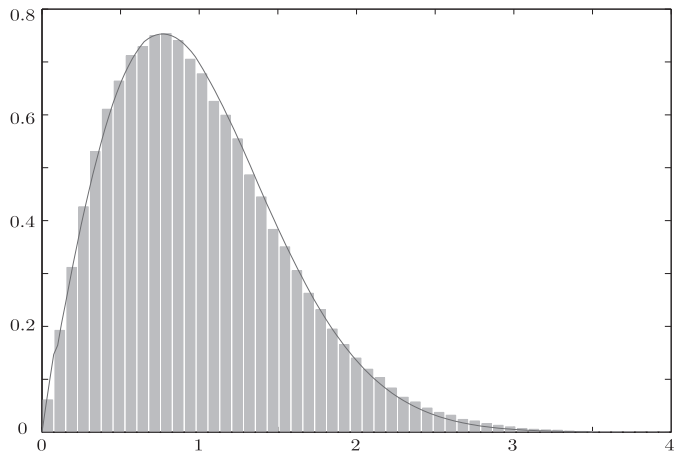}
\caption[Distribution of scaled eigenvalue spacings for a random
sample of matrices having the form $C+P^TTP$]{Distribution of scaled
eigenvalue spacings over the bulk of the spectrum for a random
sample of 300 matrices having the form $C+P^TTP$, plotted with the
GOE bulk scaling limit.  The corresponding graphs have 600 vertices
each.\label{EmpSpacGOE}}
\end{center}
\end{figure}
We generated large samples of random matrices having the form
$C+P^TTP$, using Matlab's random permutation generator.  The
limiting eigenvalue statistics of our random matrices are the same
as those of much larger, and more widely studied classes. In
particular, the limiting eigenvalue density is described by the
McKay density \cite{McKay}, and the distribution of scaled
eigenvalue spacings\footnote{Spacings are found by ordering the
eigenvalues in ascending order and finding the distance between
consecutive eigenvalues.  These distances are then scaled to have
mean one and integrated density one.} from the bulk of the spectrum
appears to be that of the Gaussian Orthogonal Ensemble (GOE). Figure
\ref{EmpMcKay} shows the empirical eigenvalue density for a random
sample of $C+P^TTP$ matrices. Figure \ref{EmpSpacGOE} shows the
scaled eigenvalue spacing distribution over the bulk of the spectrum
for the same sample.

The McKay density formula,
\begin{equation}
f_k(x)=\left\{\begin{array}{ll}
\frac{k\sqrt{4(k-1)-x^2}}{2\pi(k^2-x^2)} & \mbox{for }|x|\le2\sqrt{k-1}\\
0 & \mbox{otherwise}
\end{array}\right\}
\end{equation}
describes the limiting eigenvalue density for sequences of of
$k$-regular graphs which look like trees at most vertices.
\begin{thm}[McKay]
Given a sequence $X_1,X_2,\ldots$ of $k$-regular graphs for some
$k\ge2$, let $|X_n|$ denote the number of vertices in graph $X_n$,
and $w_r(X_n)$ the number of closed walks of length $r$.  If
$|X_n|\rightarrow\infty$ as $n\rightarrow\infty$, and
\begin{equation}
\lim_{n\rightarrow\infty}\frac{w_r(X_n)}{|X_n|}=0\qquad \forall
r\ge3
\end{equation}
then the limiting probability density function describing the
eigenvalues of $X_n$ is $f_k(x)$. \label{McKay}
\end{thm}
McKay proved that for a deterministic sequence of $k$-regular
graphs, if the limiting eigenvalue density is not described by
McKay's density formula, then the ratio of the number of closed
walks of length $r$ to the number of vertices must not approach
zero, for some value of $r$.

The Gaussian Orthogonal Ensemble is the set of $N\times N$ matrices
whose components are independent (up to the symmetry requirement)
random variables of mean zero and standard deviation one. Numerical
studies have revealed that the distributions describing the
eigenvalue statistics of the Gaussian Orthogonal Ensemble (GOE),
describe the eigenvalue statistics of most real symmetric random
matrices having independent entries.  However, the entries of our
random matrices are highly dependent.  Since these matrices
represent three-regular graphs, the entries of each row and each
column must sum to three.  Furthermore, the cycle subgraph ensures
that two of the non-zero entries of each row and each column come
from the first sub/super-diagonals of the matrix.

Given our random matrices describe one-face maps, a natural question
is whether the eigenvalue statistics depend on the genus of the
underlying map.  To investigate this question we must generate
matrices corresponding to random one-face maps having $N$ edges and
genus $g$.

\section{Constructing genus zero maps}
\label{Constructing} The number of genus $g$ one-face maps having
$N$ edges in the embedded graph is given by \cite[Theorem 2]{HarZag}

\begin{equation}\label{epsg}
\epsilon_{g}(N)=\frac{(2N)!}{(N+1)!(N-2g)!}\times\text{the
coefficient of $x^{2g}$ in }\left(\frac{x/2}{\tanh
x/2}\right)^{N+1}.
\end{equation}
The vast majority of one-face maps having $N$ edges are of genus
slightly less than $\lfloor N/2 \rfloor$.  For a given $N$, random
selections of one-face maps will only sample genera near this upper
bound.  For large $N$, the distribution of genera makes it
infeasible to generate a random sample of low genus (relative to
$\lfloor N/2 \rfloor$) by generating random one-face maps and
sorting out the maps with desired genus. We need an algorithm which
will generate random one-face maps of a specified genus.  In the
genus zero case, the bijection between one-face maps and pair
partitions provides such an algorithm.

As previously described, the important feature in the representation
of one-face maps by three-regular graphs of the form $C+P^TTP$ is
the toothpick subgraph, $P^TTP$.  This portion of the graph defines
the edge glueings of the $2N$ polygonal representation. This
information can be encoded in a pair partition diagram. If vertex
$i$ and vertex $j$ are adjacent in the toothpick subgraph, then node
$i$ and node $j$ are paired in the associated pair partition.  There
is a bijection between the set of non-crossing pair partitions and
genus zero one-face maps.  Using this bijection, we have developed
an algorithm to generate the adjacency matrices of three-regular
graphs representing random genus zero one-face maps.  The generation
of one-face maps having specified genus $g>0$ and $N$ edges remains
an open problem.

\begin{figure}[ht]
\begin{center}
\includegraphics{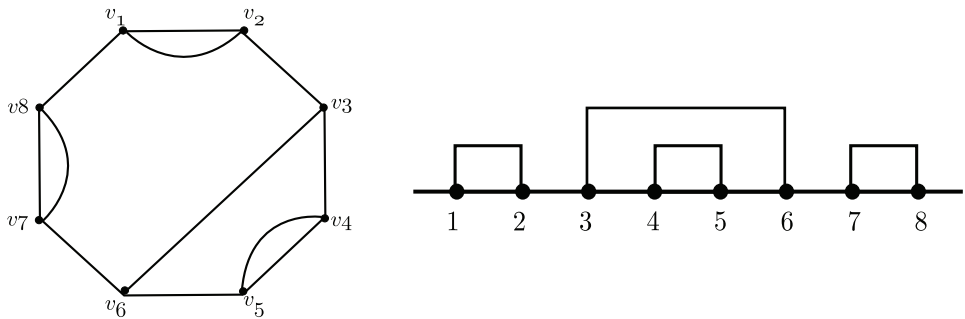}
\caption{Three-regular graph representing a genus zero one-face map
and the associated non-crossing pair partition\label{ExNCPP}}
\end{center}
\end{figure}

Our algorithm generates a random non-crossing pair partition and
then uses the node identifications of this pair partition to assign
adjacencies in the toothpick subgraph of the associated
three-regular graph.  Non-crossing pair partitions are formed
recursively.  The node paired with node one is randomly chosen
according to the probability mass function described below.  This
pairing partitions the remaining nodes into two sets, one within the
pairing and one external to it.  The first node of each of these
sets is then paired using the same algorithm.  This process is
repeated until all nodes have been paired.

\begin{figure}[ht]
\begin{center}
\includegraphics{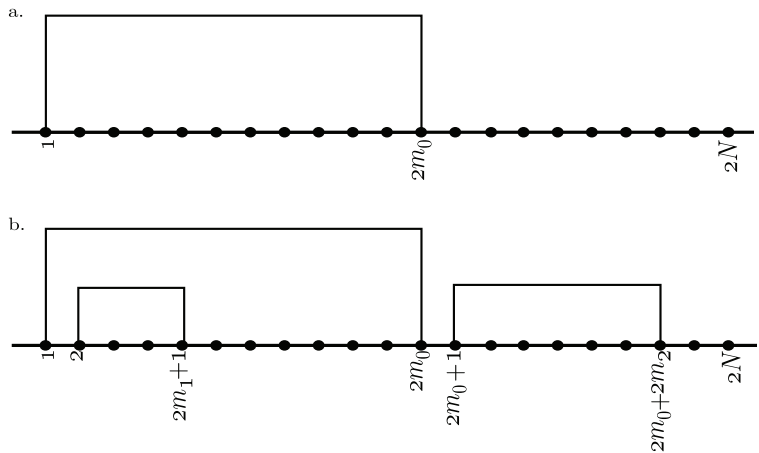}
\caption{a. The partially completed pair partition after one call to
our algorithm. The $m_0$ is randomly chosen according to the
\ensuremath{\mathit{pmf}}, and node 1 is paired with node $2m_0$. b.
The partially completed pair partition after a second round of
calls.  Elements $m_1$ and $m_2$ are randomly selected according to
the probability mass functions $pmf(m,m_0-1)$ and $pmf(m,N-m_0)$,
respectively. \label{NCPPAlg}}
\end{center}
\end{figure}

The number of non-crossing pair partition on $2N$ nodes is given by
the Catalan number, $$C_N:=\frac{1}{N+1}{2N \choose N},$$  which can
be defined recursively as \begin{equation} \label{recCN}
C_N=\frac{(4N-2)\cdot C_{N-1}}{N+1}.
\end{equation}
Let $2m$ be the random variable representing the node paired with
node one in a non-crossing pair partition. Since the pair partition
is non-crossing, odd nodes must be paired with even nodes.  The
number of non-crossing pair partitions on $2N$ nodes in which the
first node is paired with node $2m$ is equal to the number of
non-crossing pair partitions on $2m-2$ elements times the number of
non-crossing pair partitions on $2N-2m$ elements, $C_{m-1}\cdot
C_{N-m}$. Thus, the probability mass function
(\ensuremath{\mathit{pmf}}) of $m$ is a function of both $m$ and
$N$, and is given by
\begin{equation}
pmf(m,N)=\frac{C_{m-1}\cdot C_{N-m}}{C_N}.
\end{equation}
To avoid computer arithmetic errors, the $pmf(m,N)$ for $m$ greater
than five and less than $N$ is calculated using the recursive
formula derived from equation \eqref{recCN}.
\begin{equation}
pmf(m,N)=\frac{(N+1)(2N-2m-1)}{(N-m+1)(2N-1)}pmf(m,N-1).
\end{equation}
Our algorithm was implemented using Matlab.  With each recursive
call to our program, we assign a pairing to the first element of
some subset of the entire pair partition. These subsets are
determined by the previously assigned pairs. For the example shown
in Figure \ref{NCPPAlg}, after the first pairing is assigned, the
partition is divided into two subsets of lengths $2m_0-2$ and
$2N-2m_0$.  The first element of each of these subsets is paired
during the next round of calls. Since the \ensuremath{\mathit{pmf}}
used to assign these pairings is based on the subset lengths, it is
different for each recursive call.

\subsection{Proof that genus zero one-face maps are chosen with
equal probability}

As previously discussed, randomly choosing a genus zero one-face map
having $N$ edges is equivalent to randomly choosing a non-crossing
pair partition on $2N$ nodes.  The number of non-crossing pair
partitions on $2N$ nodes is $C_N$.  Thus, we want our algorithm to
construct non-crossing pair partitions on $2N$ nodes in such a way
that every partition has a $\frac{1}{C_N}$ chance of being created.

Induction Case: There is only one way to partition 2 nodes. Node 1
is paired with node 2 with probability 1, and the resulting
non-crossing pair partition is chosen with probability
$1=\frac{1}{C_1}=\frac{1}{1}$.

Suppose that for all $k$ less than $N$, the algorithm randomly
constructs non-crossing pair partitions on $2k$ nodes in such a way
that each possible partition has the equal probability of being
created, $\frac{1}{C_k}$. Let $X$ be an arbitrary non-crossing pair
partition of $2N$ elements, and let node $2m$ be the element paired
with node one in $X$. The probability our algorithm constructs
partition $X$ is equal to the probability of node 1 being paired to
node $2m$, times the probability of the particular non-crossing pair
partition $I$ on the $2m-2$ nodes internal to the pairing 1:$2m$,
times the probability of the particular non-crossing pair partition
$E$ on the $2N-2m$ nodes external to the pairing 1:$2m$. By our
supposition, our algorithm creates the internal pair partition $I$
with probability $\frac{1}{C_{m-1}}$, and the external pair
partition $E$ with probability $\frac{1}{C_{N-m}}$.  Thus, the
probability of an arbitrary non-crossing pair partition $X$ being
created is equal to
\begin{eqnarray}
P(X)&=&pmf(m,N)\frac{1}{C_{m-1}}\frac{1}{C_{N-m}}\nonumber\\
&=& \frac{C_{m-1}\cdot C_{N-m}}{C_N}\frac{1}{C_{m-1}}\frac{1}{C_{N-m}}\nonumber\\
&=&\frac{1}{C_N}
\end{eqnarray}

By induction, our algorithm randomly constructs non-crossing pair
partitions on $N$ nodes with equal probability $\frac{1}{C_N}$.

\section{Results} \label{Results}

\begin{figure}
\begin{center}
\includegraphics[width=5in]{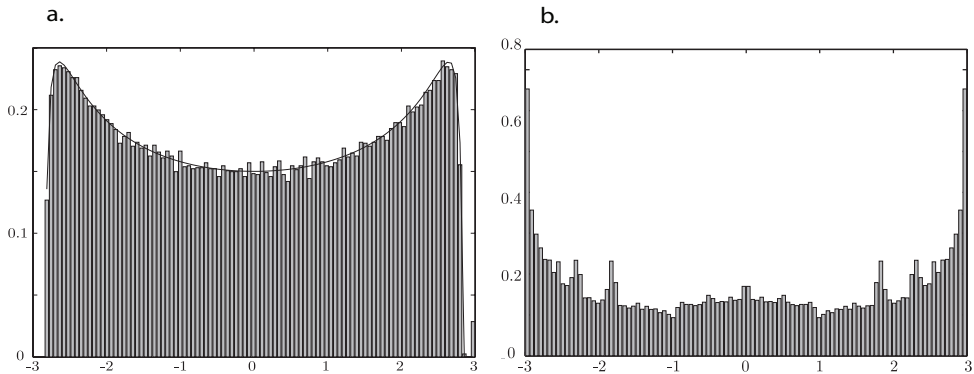}
\caption{a. The eigenvalue density of three-regular graphs
representing genus 147 one-face maps.  The sample was created by
generating random three regular graphs of the form $C+P^TTP$
corresponding to one-face maps having $N=300$ edges.  The graphs
were then sorted according to the genus of the associated map.
Thirty-six genus 147 maps were included in this sample. b. The
eigenvalue density of three-regular graphs representing genus zero
one-face maps.  The sample contains $100$ graphs, each representing
a genus zero one-face map having $N=1000$ edges.
\label{gen147gen0den}}
\end{center}
\end{figure}

\begin{figure}
\begin{center}
\includegraphics[width=5in]{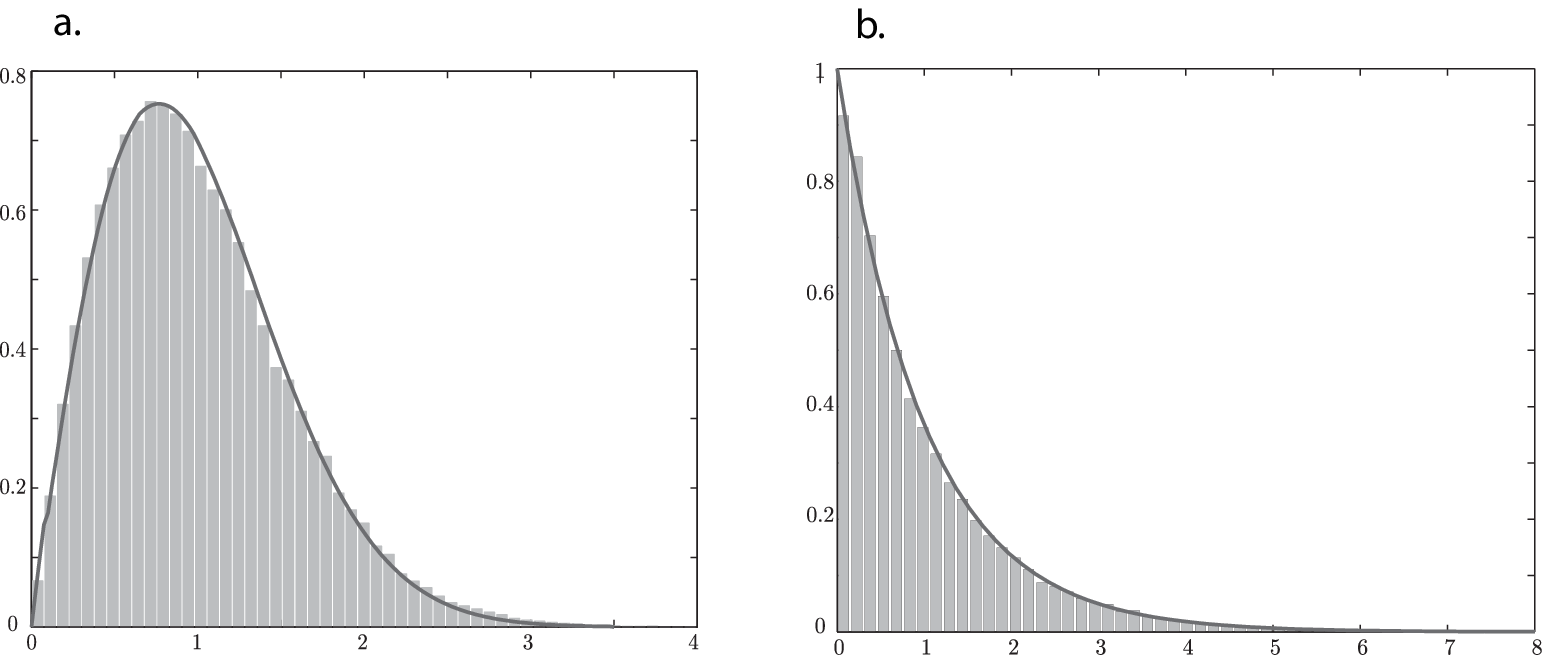}
\caption{a. The scaled eigenvalue spacing distribution of
three-regular graphs representing genus 147 one-face maps plotted
with the GOE bulk scaling limit.  The sample was created by
generating random three regular graphs of the form $C+P^TTP$
corresponding to one-face maps having $N=300$ edges.  The graphs
were then sorted according to the genus of the associated map.
One-hundred and nineteen genus 147 maps were included in this
sample. b.  The scaled eigenvalue spacing distribution of
three-regular graphs representing genus zero one-face maps, plotted
with the exponential probability density function.  The sample
contains $100$ graphs, each representing a genus zero one-face map
having $N=800$ edges. \label{gen147gen0spac}}
\end{center}
\end{figure}

We can examine the eigenvalue statistics of one-face maps having
specified genus near $\lfloor N/2\rfloor$ by forming large samples
of random one-face maps having $N$ edges and sorting by genus.
Figures \ref{gen147gen0den} and \ref{gen147gen0spac} show the
eigenvalue density and spacing distribution for a random sample of
genus 147 one-face maps generated in this way.  Examining the
eigenvalue statistics of one-face maps having genus near $\lfloor
N/2\rfloor$, we find that for large $N$ the eigenvalue density
appears to be McKay, and the scaled spacing density appears to be
that of the Gaussian Orthogonal Ensemble.  These results are
unsurprising given that samples of random one-face maps for large
$N$ exhibit the same properties and are comprised primarily of maps
having genera near $\lfloor N/2\rfloor$.


\begin{figure}
\begin{center}
\includegraphics{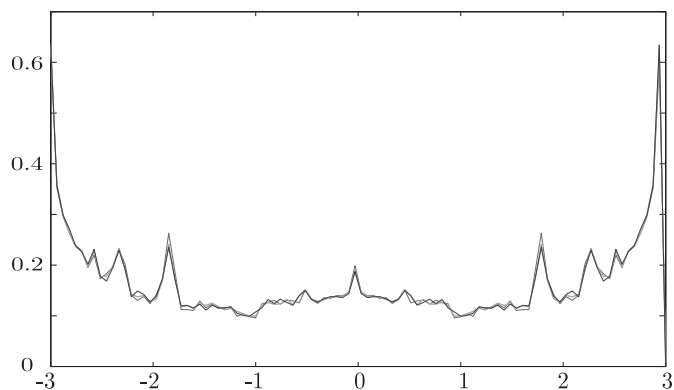}
\caption[Eigenvalue density curves for random samples of genus zero
one-face maps having $N$ edges for various values of $N$]{The
eigenvalue density curves for random samples of genus zero one-face
maps having $N=50$, $N=200$, and $N=400$ edges.  Each sample
contains 500 genus zero one-face maps. \label{gzEDcomp}}
\end{center}
\end{figure}

Using the previously described algorithm, we generated samples of
random genus zero one-face maps for large $N$, and examined the
eigenvalue statistics. The results are strikingly different from
those found for one-face maps having genus near $\lfloor
N/2\rfloor$.  The distribution of eigenvalues for a random sample of
genus zero one-face maps is plotted in Figure \ref{gen147gen0den}.
In Figure \ref{gzEDcomp} the eigenvalue density curves for various
samples of genus zero one-face maps are plotted on the same axis.
Although the number of embedded edges $N$ is different for each
sample, the eigenvalue densities are nearly identical.  While this
distribution is the same for all large values of $N$, it does not
appear to be described by any established density formula. Examining
the eigenvalues themselves, we find that the peaks observed in the
density curve near $0, \pm 0.5, \pm 1.8$ and $\pm 2.3$ are not the
result of repeated eigenvalues, but rather a decrease in spacing
between consecutive eigenvalues.  The non-crossing structure of the
toothpick portion of the three-regular graphical representation of
genus zero one-face maps ensures that the graphs are bipartite. The
symmetry in the distribution of eigenvalues about $\lambda =0$ is
due to this bipartite structure.

\begin{figure}
\begin{center}
\includegraphics{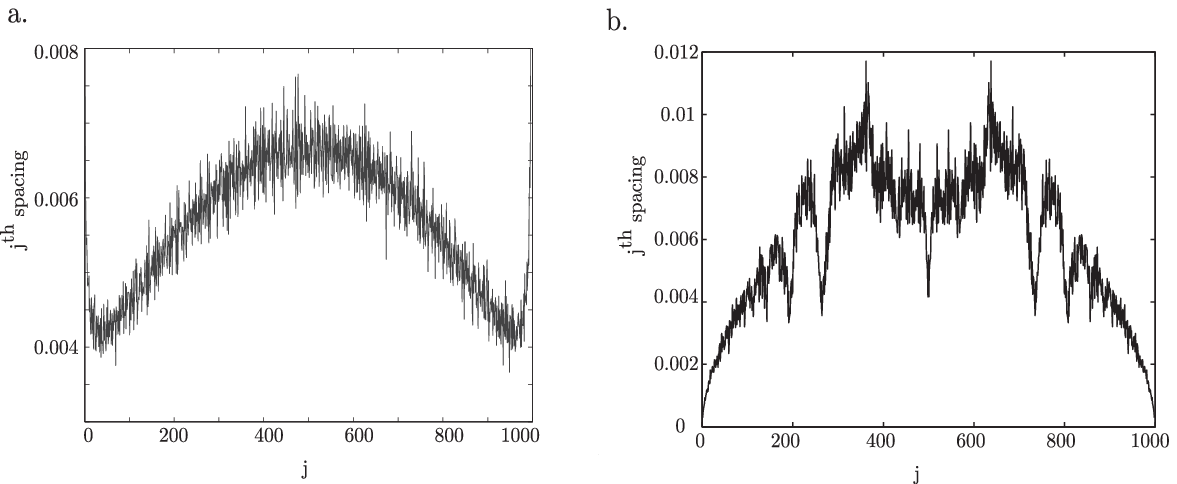}
\caption{a. The mean $j^{th}$ eigenvalue spacing as a function of
$j$ for a sample of 100 three-regular graphs representing one-face
maps of genera near $g=250$.  b. The mean $j^{th}$ eigenvalue
spacing as a function of $j$ for a sample of $100$ three-regular
graphs representing genus zero one-face maps.  Both samples contain
graphs having $1000$ vertices.  \label{meanjthSpac}}
\end{center}
\end{figure}

The difference in the distribution of eigenvalues for the genus zero
and genus $\sim\lfloor N/2\rfloor$ cases is also seen in the
corresponding mean $j^{th}$ spacing curves.  Figure
\ref{meanjthSpac} shows the mean $j^{th}$ eigenvalue spacing as a
function of $j$ for a sample of three-regular graphs representing
genus zero one-face maps, and a sample of three-regular graphs
representing one-face maps of genera near
$\lfloor\frac{N}{2}\rfloor$.

\begin{figure}
\begin{center}
\includegraphics{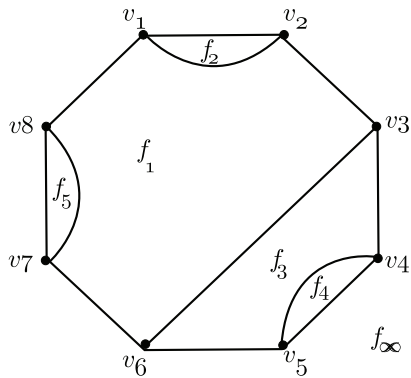}
\caption[Three-regular graphical representation of a genus zero
one-face map]{The three-regular graphical representation of a genus
zero on-face map with its faces labeled.\label{Genus0Ex}}
\end{center}
\end{figure}

It is clear from these results that, in the genus zero case, the
distribution of eigenvalues is not described by the McKay density
formula.  In fact, three-regular graphs representing genus zero one
face maps violate the criteria of McKay's theorem (\ref{McKay}).
These graphs are always planar, i.e. they can be drawn in the plane
in such a way that the edges of the graph do not intersect except at
vertex points.  A planar drawing of a graph divides the set of
points of the plane not lying on the graph into regions called {\it
faces}. Figure \ref{Genus0Ex} shows the three-regular graph
representing a genus zero one-face map with its faces labeled. One
of the faces will be unbounded.  This face is called the {\it
infinite} face. All other faces are bounded by the edges of the
graph, and are called {\it internal} faces.

For genus zero one-face maps, the embedded graph is a planar tree,
and each internal face of the associated three-regular graph
corresponds to a vertex of this embedded tree. The degree of this
embedded vertex is equal to $\frac{1}{2}$ the face boundary length.
Thus, the distribution of vertex degrees in the genus zero embedded
trees determines the distribution of face boundary lengths in the
associated three-regular graph. Furthermore, since each face
boundary forms a cycle in the graph, the number of vertices of
degree $d$ in the tree provides a lower bound for the number of $2d$
cycles in the three-regular graph.

Let $\mathcal{X}_N$ denote the set of three-regular graphs having
$2N$ vertices which represent genus zero one-face maps. Drmota and
Gittenberger \cite{Drmota} proved that the number of degree $k$
vertices in a random planar tree having $N$ edges is asymptotically
$\frac{N}{2^k}$. Thus, the expected number of $2k$-cycles in a
random element of $\mathcal{X}_N$ as $N\rightarrow\infty$ is
$$\lim_{N\rightarrow\infty}\mathbb{E}(c_{2k})\ge\frac{N}{2^k}.$$
Since the number of closed walks of length $2k$ is greater than or
equal to the number of cycles of length $2k$, the set of
three-regular graphs representing genus zero one-face maps contains
a subset $\{X_1,X_2,\ldots\}$ of nonzero measure for which $$
\lim_{N\rightarrow\infty}\frac{w_{2k}(X_N)}{|X_N|}\ge\lim_{N\rightarrow\infty}\frac{c_{2k}(X_N)}{2N}>
0.$$ This violates the criteria of McKay's Theorem.

The distribution of scaled eigenvalue spacings for a sample of
random genus zero one-face maps having $N=800$ edges is plotted in
Figure \ref{gen147gen0spac}. The scaled spacings in the genus zero
case are not described by the GOE density formula which was shown
numerically to describe the bulk scaling limit for genera near
$\lfloor N/2\rfloor$, but rather appear to be described by the
exponential probability density function.

\section{Conclusion and conjecture}\label{VarMean}

The numerical results suggest that the eigenvalue statistics depend
on the rigidity of the associated three-regular graph. For one-face
maps having small genus relative to $N$, a great deal of structure
is imposed on the corresponding three-regular graph.  We believe the
striking difference between the eigenvalue statistics of genus zero
maps and genus $\sim\lfloor N/2\rfloor$ is a result of this imposed
structure. For any fixed genus, as the size of the map increases
(i.e. the number of vertices in the associated three-regular graph
grows), more structure will be imposed on the associated
three-regular graph.  This leads to the conjecture that for fixed
genus $g$, in the limit as $N\rightarrow\infty$, the scaled spacing
distribution will be exponential.  As we have seen, for $N\sim 2g$,
the scaled spacing distribution appears to be GOE.  Thus the
distribution must undergo a shift from GOE to exponential as $N$
increases.  Similar shifts in probability density functions have
been observed in biological networks and quantum chaotic systems
\cite{Bog,Luo}.

%


\end{document}